\documentclass{amsart} 
\newtheorem{theorem}{Theorem}[section]
\newtheorem{lemma}[theorem]{Lemma}
\newtheorem{corollary}[theorem]{Corollary}
\newtheorem{proposition}[theorem]{Proposition}
\newtheorem{remark}[theorem]{Remark}

\theoremstyle{definition}
\newtheorem{definition}[theorem]{Definition}
\newtheorem{example}[theorem]{Example}
\newtheorem{remark/example}[theorem]{Remark/Example}

\newtheorem{question}[theorem]{Question}

\newtheorem{claim}[theorem]{Claim} 

\let\oldlabel=\label
\def\prellabel{\marginparsep=1em\marginparwidth=44pt
 \def\label##1{\oldlabel{##1}\ifmmode\else\ifinner\else
 \marginpar{{\footnotesize\ \\ \tt
 ##1}}\fi\fi}}

\numberwithin{equation}{section}

\newcommand{\NN}{ {\bf N} }
\newcommand\ZZ{ {\bf Z} }

\newcommand\KK{ {\bf K} }

\newcommand\M{{\bf m}}

\newcommand\C{\mathcal C}
\newcommand\D{{\mathcal C}^*}
\newcommand\G{\mathcal G}
\newcommand\A{\mathcal A}
\newcommand\Q{\mathcal Q}

\newcommand{\Rees}{\operatorname{Rees}}

\newcommand{\GCD}{\operatorname{GCD}}

\newcommand{\length}{\operatorname{length}}

\newcommand{\Image}{\operatorname{Image}}
\newcommand{\reg}{\operatorname{reg}}
\newcommand{\gr}{\operatorname{gr}}
\newcommand{\Ass}{\operatorname{Ass}}
\newcommand{\HF}{\operatorname{HF}}
\newcommand{\HS}{\operatorname{HS}}
\newcommand{\Min}{\operatorname{Min}}
\newcommand{\sat}{\operatorname{sat}}
\newcommand{\comp}{\operatorname{comp}}

\numberwithin{equation}{section}

\begin{document}
\title{Integrally closed and componentwise linear ideals}
\author{ A. Conca, E. De Negri, M. E. Rossi }
\address{Dipartimento di Matematica, Universit\'a di Genova, Via Dodecaneso 35, I-16146 Genova, Italy }
\email{conca@dima.unige.it, denegri@dima.unige.it, rossim@dima.unige.it }
\subjclass[2000]{Primary 13B22; Secondary 13D02}
\date{}
\thanks{Research partially supported  by the ``Ministero dell'Universit\'a e della  Ricerca Scientifica" in the framework of the National Research Network (PRIN 2005)  ``Algebra commutativa, combinatorica e computazionale"}

\keywords{integrally closed ideals, associated graded rings, componentwise linear ideals, factorization properties}

\begin{abstract}
In a two dimensional regular local ring  integrally closed ideals have a unique factorization property and their associated graded ring is  Cohen-Macaulay. In higher dimension these properties do not hold and the goal of the paper is to identify a subclass of integrally closed ideals for which they do. We restrict our attention to $0$-dimensional homogeneous ideals in polynomial rings $R$ of arbitrary dimension. We identify a class of integrally closed ideals,  the Goto-class $\G^*$, which  is closed under product and it has a suitable unique factorization property. Ideals in $\G^*$ have a Cohen-Macaulay associated graded ring if either they are monomial or $\dim R\leq 3$. Our approach is based on the study of the relationship between the notions of integrally closed, contracted, full and componentwise linear ideals. 
\end{abstract}

\maketitle
\section{Introduction}
\label{intro}
Thanks to the work of Zariski, integrally closed ideals of two-dimensional regular local rings $(R,\M)$ are well-understood. In such rings the product of integrally closed ideals is integrally closed and there is a unique factorization property for integrally closed ideals into product of simple integrally closed ideals. In higher dimension, these properties no longer hold, see the examples in \cite{C,C3,H,L}. 
The identification of  analogues of Zariski's results is an active research area. In this direction we mention the work of Cutkosky \cite{C,C1,C2,C3}, Deligne \cite{D}, Huneke \cite{H} and Lipman \cite{L}. Several authors considered other  related problems, as for instance the description of integrally closed ideals $I$ such that $I\M$ is integrally closed as well, see \cite{CGPU,DC1,DC2,EM,Ga,Ga1,HH}. 

 In this paper we deal with homogeneous ideals of $R=K[x_1,\dots,x_n]$, the polynomial ring over a field $K$.  For an ideal $I$ we denote by $o(I)$ the order or initial degree of $I$, by $I_j$ the homogeneous component of degree $j$ of $I$ and by $I_{\langle j\rangle}$ the ideal generated by $I_j$. We set  $\M=(x_1,\dots,x_n)$.

Our goal is to identify a class of $\M$-primary integrally closed ideals of $R$ which behaves, as much as possible, as the class of integrally closed ideals in dimension $2$. To this end, we study the relations between four properties of ideals: 1) being integrally closed, 2) being componentwise linear, 3) being contracted (from a quadratic extension), 4) being $\M$-full. 
It turns out that 1) implies 3), that 2) implies 3) and that 3) implies 4). 
Also, for ideals $I$ such that $I+(\ell)=\M^{o(I)}+(\ell)$ for some linear form $\ell$, one has that 4) implies 2). 
 
We then consider the class $\C$ of the $\M$-primary ideals of $R$ satisfying $I+(\ell)=\M^{o(I)}+(\ell)$ for some linear form $\ell$ and having property 4), (equivalently 3) or 2)). Denote by $\D$ the set of the ideals in $\C$ that are integrally closed.
 We prove that $\C$ is closed under product and integral closure, see Proposition \ref{prodotto}. Further, we prove in Theorem \ref{factorization} that $\C$ has a factorization property that looks like Zariski's factorization for contracted ideals in dimension $2$ \cite[Appendix 5, Thm.1]{ZS}. An important role in Zariski's factorization theorem is played by the characteristic form $g(I)$ defined has the $\GCD$ of the forms of degree $o(I)$ in $I$. Given $I\in \C$ for every $j\in \NN$ we define $Q_j(I)$ to be the saturation of $I_{\langle j+o(I)\rangle}$. 
 In our context, the characteristic form is replaced by the ideal $Q_0(I)$. 
 
 We show that given $I\in \C$, one has $I\in \D$ iff $I\M\in \D$. But, unfortunately, $\D$ is not closed under product. 
We then consider the Goto-class $\G$ defined as the set of the ideals $I\in\C$ such that for every $j$ the primary components of $Q_j(I)$ are powers of (necessarily $1$-dimensional) geometrically prime ideals. 
Integrally closed complete intersections, characterized by Goto \cite{G}, are in $\G$, see Theorem \ref{easy}. 
We prove in Proposition \ref{1/2G} that $\G$ is closed under product and that it is compatible with the factorization of $\C$. 
We define $\G^*$ to be the set of the integrally closed ideals of $\G$. We then show that $\G^*$ is closed under product and has a unique factorization property, see Theorem \ref{mainG}. The simple elements in $\G^*$ have a ``simple" description: up to a change of coordinates, they are of the form $\overline{(x_1^d,\dots, x_{n-1}^d, x_n^t)}$ for coprime $d,t$ with $d<t$. 
Lipman and Teissier \cite{LT} and Huneke \cite{H2} proved that integrally closed ideals in two dimensional regular local rings have a Cohen-Macaulay associated graded ring. It is natural to ask whether the same holds for ideals of $\G^*$. We conclude the paper by showing that if $I\in \G^*$ and either $I$ is monomial (e.g.~$Q_0(I)$ has at most two minimal primes) or $\dim R\leq 3$, then the associated graded ring $\gr_I(R)$ is Cohen-Macaulay, see Corollary \ref{megliocheniente} and Theorem \ref{3CM}. 

\section{ $\M-$full, contracted and componentwise linear ideals}
\label{M-full}
 
 Throughout the paper let 
$R=K[x_1,\dots,x_n]$ be a polynomial ring over a field $K$, and $\M=(x_1,\dots,x_n)$. 
All the ideals we deal with are homogeneous (with few exceptions).

Let $I$ be an ideal of $R$. Denote by $\mu(I)$ the minimum 
number of generators of $I$ and by
$o(I)$ the initial degree (or the order) of $I$, that is the least degree of non-zero elements in $I$.

 In this section we discuss the relations between $\M$-full, contracted and componentwise linear ideals. First we introduce some notation and recall definitions.
Denote by $\beta_{ij}(I) $ the $ij-$th graded Betti number of $I$ as an $R$-module. 
The Castelnuovo-Mumford regularity of $I$ is given by 
$$ \reg(I) = \max \{ j-i : \beta_{ij}(I)\neq 0 \}.$$ 
 The ideal $I$ has a linear resolution if $\reg(I)=o(I)$. For general facts on the Castelnuovo-Mumford regularity and its characterization in terms of local cohomology we refer the reader to \cite{E}. 
For every integer $j$ denote by $I_j$ the $K$-vector space of the forms of degree $j$ in $I$, and by $ I_{\langle j\rangle}$ the ideal generated by the elements of $I_j$. 
 The ideal $I_{\langle j\rangle}$ has a linear resolution for $j\geq \reg(I)$.
 
 Given two ideals $I$ and $J$, we set $I:J^{\infty}=\cup_kI:J^k$. We denote by $I^{\sat}$ the saturation of $I$ with respect to $\M$, that is $$I^{\sat}=I:\M^{\infty}.$$
  For short we will denote the ideal  $ ( I_{\langle j\rangle})^{\sat}$ by 
 $ I_{\langle j\rangle}^{\sat} $.

\begin{definition}
\label{defcl} 
 An ideal $I\subset R$ is said to be componentwise linear if $ I_{\langle d\rangle}$ has a linear resolution for every $d\in \NN$. 
\end{definition}

For every non-zero linear form $\ell$ in $R$ we consider the quadratic transform $S$ of $R$ associated to $\ell$. By definition $S=R[\M/\ell]=\cup_{k\in \NN} \M^k/\ell^k$. 

 \begin{definition} 
\label{defcontratti} An ideal $I\subset R$ is said to be contracted (from a quadratic extension) if there exists a non-zero linear form $\ell$ in $R$ such that $I=IS\cap R$, where $S=R[\M/\ell]$. 
\end{definition}

\begin{proposition} 
\label{ells}
Let $\ell$ be a non-zero linear form in $R$ and $I\subset R$ an ideal. Set $S=R[\M/\ell]$ and $J= IS\cap R $. We have: 
\begin{itemize}
\item[(1)] $J=\cup_{k\in \NN} ( I\M^k:\ell^k )$. 
\item[(2)] $J$ is homogeneous. 
\item[(3)] $J_j=(I_{\langle j\rangle}^{\sat}:\ell^{\infty})_j$. 
\end{itemize}
 \end{proposition}

\begin{proof} (1) follows immediately from the fact that $IS=\cup_k I\M^k/\ell^k$. Then (2) follows from (1). To prove (3) consider $f\in R$ homogeneous of degree $j$. We have $f\in J_j$ iff $f\ell^k\in (I\M^k)_{j+k}$ for every $k\gg 0$. 
Since $(I\M^k)_{j+k}=(I_{\langle j\rangle})_{j+k}$ we have $f\in J_j$ iff $f\ell^k\in I_{\langle j\rangle}$ for every $k\gg 0$. Hence $f\in J_j$ iff $f\in I_{\langle j\rangle}:\ell^{\infty}= I_{\langle j\rangle}^{\sat}:\ell^{\infty}$. 
\end{proof}

In the following we denote by $\Ass(M)$ the set of the associated prime ideals of an $R$-module $M$. 
\begin{definition} Let $I$ be an ideal of $R$. We set 
$$\Ass^{\comp}(R/I)=\cup_{j\geq o(I)} \Ass(R/I_{\langle j\rangle}).$$
\end{definition} 

\begin{lemma} Let $I$ be an ideal of $R$ with generators in degrees $d_1, \dots,d_p$, $d_1<\dots<d_p$. We have 
$$\Ass^{\comp}(R/I)=\Ass(R/I_{\langle d_1 \rangle})\cup \dots \cup \Ass(R/I_{\langle d_p \rangle}) \cup \{\M\}.$$
In particular, $\Ass^{\comp}(R/I)$ is finite. 
\end{lemma} 

\begin{proof} The assertion follows immediately by observing that if $I$ has no generators in degree $j+1$, then 
$I_{\langle j+1\rangle}=I_{\langle j \rangle} \cap \M^{j+1}$. 
\end{proof} 

\begin{definition} Let $I$ be an ideal. We denote by $U(I)$ the (finite) union of the prime ideals in $\Ass^{\comp}(R/I)\setminus\{\M\}$.
\end{definition}
 
\begin{proposition} \label{dk} Let $I$ be an ideal with generators in degrees $d_1, \dots,d_p$ with $d_1<\dots<d_p$ and set $d_{p+1}=\infty$. The following conditions are equivalent: 
\begin{itemize}
\item[(1)] $I$ is contracted from $R[\M/\ell]$ for some non-zero linear form $\ell$. 
\item[(2)] $I$ is contracted from $R[\M/\ell]$ for every non-zero linear form $\ell$ with $\ell \not\in U(I)$. 
\item[(3)] $( I_{\langle j\rangle}^{\sat} )_j =I_j$ for every $j\in \NN$. 
\item[(4)] $( I_{\langle d_k\rangle}^{\sat} )_j =I_j $ for every $j$ with $d_k \le j <d_{k+1}$ and $k=1, \dots, p$. 
\end{itemize} 
\end{proposition}
\begin{proof} Obviously  (2) implies (1). That (1) implies (3) follows from $I_j=(I_{\langle j\rangle}^{\sat}:\ell^\infty)_j$, which holds by \ref{ells}, and $(I_{\langle j\rangle}^{\sat}:\ell^\infty)_j \supseteq (I_{\langle j\rangle}^{\sat})_j \supseteq I_j$. For (3) implies (2) one notes that if $\ell\not\in U(I), $ then we have $I_{\langle j\rangle}^{\sat}:\ell^\infty=I_{\langle j\rangle}^{\sat}$ and by assumption $( I_{\langle j\rangle}^{\sat} )_j =I_j$. It follows then from \ref{ells} that $I$ is contracted from $R[\M/\ell]$. Finally, that (3) and (4) are equivalent follows from the observation that if $I$ has no generators in degree $j+1, $ then $I_{\langle j+1\rangle}=I_{\langle j\rangle}\cap \M^{j+1}$ and hence $I_{\langle j+1\rangle}^{\sat}=I_{\langle j\rangle}^{\sat}$.
\end{proof}

\begin{proposition} \label{clToco} Every componentwise linear ideal of $R$ is contracted.
\end{proposition}
\begin{proof} Since $I$ is componentwise linear, we have $\reg(I_{\langle j\rangle})=j$ for every $j$ and hence $I_j= (I_{\langle j\rangle})_j = (I_{\langle j\rangle}^{\sat})_j$. The result follows by \ref{dk} (2). \end{proof}
 
 In dimension $3$ or higher contracted ideals need not be componentwise linear. 
 
 \begin{example}\label{exe1}
 $(x_1^2,x_2^2)$ is contracted but not componentwise linear in $K[x_1,x_2,x_3]$. 
 \end{example} 
 
The following definition is due to Rees. We adapt it to the graded case. 

\begin{definition} 
\label{deffull} 
An ideal $I\subset R$ is said to be $\M-$full if there exists 
 a non-zero linear form $\ell$ in $R$ such that $I\M : \ell = I $. 
\end{definition}

Ideals which are $\M$-full are studied in \cite{W1,W2,W3,G}.
It is easy to see that if $I$ is $\M-$full, then $I : \ell = I : \M$. Moreover, if $I$ is $\M$-full then $I\M : \ell = I$ holds for a general linear form $\ell$. By \ref{ells} we have immediately that: 
 
 \begin{proposition}\label{coTofu}
 \label{conMfull} 
 Every contracted ideal of $R$ is $\M$-full. 
\end{proposition}

The following example shows that the converse of \ref{conMfull} does not hold. 

 \begin{example} \label{exe2} The ideal $I=(x_1^3, x_2^3,x_1^2x_3)+(x_1,x_2,x_3)^4$ of $K[x_1,x_2,x_3]$ is $\M$-full. But $I$ is not contracted and $I \M$ is not $\M-$full. 
 \end{example}

We recall that an element $a$ of $R$ is said to be integral over $I$ if it satisfies an equation of the form
$a^t +r_1a^{t-1} + \cdots +r_t = 0$, with $r_i\in I^i$ for every $i=1,\ldots,t$. The elements of $R$ which are integral over $ I$ form an ideal, the integral closure of $I$, denoted by $\overline I$. An ideal is said to be integrally closed if it coincides with its integral closure.

\begin{proposition} \label{ic} Let $\ell\in R_1\setminus U(I)$ and $S=R[\M/\ell]$. Then
 $$ I \subseteq I S \cap R \subseteq \overline I. $$
\end{proposition}
\begin{proof} By \ref{ells} we have for every $j$
$$ (IS \cap R)_j = ( I_{\langle j\rangle}^{\sat} )_j.$$
Hence for every $f \in (IS \cap R)_j $ we have $f \M^k \subseteq I_{\langle j\rangle} \M^k$ for some $k$. The ``determinant trick" implies that $f \in \overline{ I_{\langle j\rangle}}$. In particular, $f \in \overline I$. 
 \end{proof}
 As a corollary we have: 
 \begin{corollary} \label{iic} 
Every integrally closed ideal of $R$ is contracted. 
\end{corollary}

Under the assumption that $I$ is $\M$-primary \ref{iic} is proved in \cite[Lemma 3.3]{DC1}. Further in \cite[2.4]{G} it is proved that integrally closed ideals are $\M$-full in a much more general context. 
 Summing up, we have seen that the following implications hold: 

$${\rm{ Componentwise\ linear }} \Longrightarrow {\rm{ Contracted \ }} \Longrightarrow\ \ \M-{\rm{ full \ }}$$
$$ \ \ \ \ \ \ \ \ \ \ \ \ \ \ \ \ \ \ \Uparrow $$
$$ \ \ \ \ \ \ \ \ \ \ \ \ \ \ \ {\rm{ Integrally\ \ closed\ }} $$
 
In dimension $2$, componentwise linear, contracted and $\M$-full are equivalent properties, but, as seen in \ref{exe1} and \ref{exe2}, in dimension $3$ and higher they differ. 

For an $R$-module $M$ we denote by $\length(M)$ its length. 

\begin{lemma}\label{mu}
Let $I$ be an $\M$-primary, $\M$-full ideal of order $d$. For every ideal $J$ containing $I$ and for every $\ell $ such that $ I \M : \ell =I$ one has
$$\mu(I) - \mu(J) = \length( \M J / \M I + \ell J).$$
It follows that $\mu(I)\geq \mu(J)$ and, in particular, $\mu(I)\geq \mu(\M^d)$. 
\end{lemma}
\begin{proof} See \cite[Lemma 2.2. (2)]{G}.\end{proof} 

One says that $I$ has the Rees property if $\mu(I)\geq \mu(J)$ for every ideal $J\supseteq I$. Under the assumption that $I$ is componentwise linear ideal, the inequality $\mu(I)\geq \mu(\M^d)$ is proved in \cite[3.4]{CHH}. A sort of Rees property is still valid for $\M$-full ideals not necessarily $\M$-primary. We refer to \cite[3.2]{CHH} for the corresponding result for componentwise linear ideals.

\begin{proposition} 
\label{generalfull}
Let $I $ and $J$ be ideals of $R$. Assume that $I$ is $\M$-full, $I \subseteq J$ and $I_t=J_t $ for $t \gg 0$. Then $\mu(I) \ge \mu(J)$. 
\end{proposition}

\begin{proof} 
First we remark that if $I$ is $\M$-full, then $I+\M^t$ is $\M$-full for every integer $t > 0$. Now, since $I + \M^t\subseteq J+ \M^t$ and $I+\M^t$ is $\M$-primary and $\M$-full ideal, it follows that $\mu(I + \M^t) \ge \mu(J + \M^t) $ by \ref{mu}. Since $I_t=J_t$ for $t\gg 0$, the inequality $\mu(I + \M^t) \ge \mu(J + \M^t) $ for $t \gg 0$ implies that $\mu(I) \ge \mu(J)$. 
\end{proof}

\begin{proposition} 
\label{elle}
 Let $I\subset R$ be an ideal of order $d $ and let $\ell $ be a non-zero linear form. Assume that $I+(\ell) = \M^d + (\ell)$. Then 
 \begin{itemize} 
 \item[(1)] if $I$ is $\M$-primary, then $\mu(I) \le \mu(\M^d)$. 
 \item[(2)] $I=I_{\langle d\rangle}+\ell(I:\ell)$.
 \item[(3)] $\dim R/ I_{\langle d\rangle} \le 1$.
 \end{itemize} 
\end{proposition}
\begin{proof} 
(1) If $I+(\ell) = \M^d + (\ell)$ holds for a linear form, then it holds for a generic linear form. Thus we may  consider  a sequence $y_1,\dots, y_n$  of generic linear forms in $R$ with $I+(y_1) = \M^d + (y_1)$, and set 
$$\alpha_i(I)= \length([I+ (y_1,\dots, y_i) ] : y_{i+1} /[ I+(y_1,\dots, y_i)] ).$$ 
By \cite[1.2]{CHH}, we have $\mu(I) \le \sum_{i=0}^{n-1} \alpha_i(I)$. We remark that $\alpha_0(I) = \length(I:y_1/I)$. By the exact sequence:
$$0\to\frac{I:y_1}{ I}\to \frac{R}{I}\to \frac{R}{I}\to\frac{R}{ I+(y_1)}\to 0$$
it follows that $\length(I:y_1/ I)=\length(R/(I+(y_1)))$. Since $I+(y_1) = \M^d + (y_1)$, we have $\alpha_0(I) =\length(R/ \M^d + (y_1))$. Moreover for every integer $i\ge 1 $ we have 
$[(y_1,\dots, y_i) + I ] : y_{i+1} /[ (y_1,\dots, y_i) + I] = [(y_1,\dots, y_i) + \M^{d} ] : y_{i+1} /[ (y_1,\dots, y_i) + \M^{d}]$. Then $ \alpha_i(I)= \alpha_i (\M^{d}) $ and the result follows since $ \sum_{i=0}^{n-1} \alpha_i(I) = \sum_{i=0}^{n-1} \alpha_i(\M^d) $ and $ \sum_{i=0}^{n-1} \alpha_i(\M^d)= \mu(\M^{d})$. 

(2) The inclusion $\supseteq$ is obvious. To prove the other inclusion we note that by assumption $\M^d\subseteq I_{\langle d\rangle}+(\ell)$. Thus $I\subseteq \M^d+(\ell)\subseteq I_{\langle d\rangle}+(\ell)$, in particular $I \subseteq I_{\langle d\rangle}+(\ell)\cap I=I_{\langle d\rangle}+\ell(I:\ell)$. 

(3) By assumption, $\M^d\subseteq I_{\langle d\rangle}+(\ell)$, that is $\M=\sqrt{\ell} \mod I_{\langle d\rangle}$. The conclusion follows by Krull hauptidealsatz. 
\end{proof}

We are ready to prove the following theorem. 
 \begin{theorem} 
\label{teoremone} 
 Let $I$ be an $\M$-primary ideal of order $d $ such that $I+(\ell)=\M^d+(\ell)$ for some 
non-zero linear form $\ell$. The following conditions are
equivalent: 
\begin{itemize} 
 \item[(1)] $\mu(I)=\mu(\M^d)$, 
 \item[(2)] $I$ is $\M$-full, 
\item[(3)] $I$ is contracted, 
 \item[(4)] $I$ is componentwise linear.
 \end{itemize}
\end{theorem}
\begin{proof} The implications $ (4) \Longrightarrow (3) \Longrightarrow\ (2)$ hold in general by \ref{clToco}, \ref{coTofu}. 
That $(2) $ implies $(1) $ follows by \ref{elle}(1) and \ref{mu}. It remains to prove $(1) $ implies $(4)$. We may assume that $I+(\ell)=\M^d+(\ell)$ for a general linear form. With the notation of the proof of \ref{elle}, one sees that the assumption $(1)$ can be stated as $ \mu(I) = \sum_{i=0}^{n-1} \alpha_i(I)$. Then by \cite[2.3, 1.5] {CHH}, we conclude that $I$ is componentwise linear.
\end{proof}

In dimension $2$ products of contracted ideals are contracted. This is not true in higher dimension. 

\begin{example} Let $R=K[x_1, x_2,x_3]$, and $I=(x_1^2,x_1x_2^2, x_2^2 x_3^2)$. 
The ideal $I$ is componentwise linear and hence contracted and $\M$-full. But $I^2$ is not $\M$-full (therefore not contracted and not componentwise linear). 
Take $J=I+\M^5$ to get an $\M$-primary example. 
\end{example}

The following result will be useful in the next section. 

\begin{theorem} \label{Id}
Let $I,J$ be componentwise linear ideals. Let $d$ be the order of $I$ and assume that $\dim R/I_{\langle d\rangle} \le 1$. Then $IJ$ is componentwise linear.
\end{theorem}
\begin{proof} 
 First assume that $I$ is generated in degree $d$. One has $(IJ)_{d+s}=I_dJ_s $ for every $s \in \NN$. Now since $\dim R/I_{\langle d\rangle} \le 1$, by \cite[2.5]{CH}, $\reg(I_{\langle d\rangle } J_{\langle s\rangle })=d+s$. Hence $IJ$ is componentwise linear.

Assume now that $I$ has generators in various degrees. Let $y_1,\dots,y_n $ be a generic sequence of linear forms. For $1\le p \le n $ denote by $ H_1(y_1,\dots, y_p,R/IJ)$ the first homology of the Koszul complex of $R/IJ$ with respect to $y_1,\dots,y_p$. 
In order to prove that $IJ$ is componentwise linear, 
by \cite[1.5, 2.2]{CHH}, it suffices to prove that $\M H_1(y_1,\dots,y_p,R/IJ)=0$ for every $p$. Since $\dim R/I_{\langle d\rangle} \le 1 $ and $\reg(I_{\langle d\rangle} + (y_1)) \le \reg(I_{\langle d\rangle})= d$ we deduce that $I + (y_1)= \M^d + (y_1)$. Consider the Koszul complex:
$$ \KK: \ \ \ \ \ \ \ \ \ \ \ \ \dots\to R^{\binom{p}{2}}\overset{\varphi_2}\to R^p\overset{\varphi_1}\to R.$$ 
We have to prove that $\M(\alpha_1,\dots, \alpha_p)\in \Image(\varphi_2)+IJR^p$ for every $(\alpha_1,\dots, \alpha_p)\in R^p$ satisfying $\varphi_1(\alpha_1,\dots, \alpha_p)\in IJ$.

 Since $I+(y_1)=\M^d+(y_1)$, then by \ref{elle}(2), we have $I=I_{\langle d\rangle}+y_1(I:y_1)$. Thus 
 $$IJ=[I_{\langle d\rangle}+y_1(I:y_1)]J=I_{\langle d\rangle}J+y_1(I:y_1)J.$$
As consequence we may write $\alpha_1y_1+\alpha_2y_2+\cdots+\alpha_py_p=a+by_1$ with $a\in I_{\langle d\rangle}J$ and $b\in (I:y_1)J, $ that is $(\alpha_1-b)y_1+\alpha_2y_2+\cdots+\alpha_py_p\in I_{\langle d\rangle}J$ which is componentwise linear by the first part of the proof, thus $\M(\alpha_1-b,\alpha_2\dots, \alpha_p)\in \Image(\varphi_2)+I_{\langle d\rangle}JR^p\subseteq \Image(\varphi_2)+IJR^p$. The conclusion follows by noting that $\M b\in J\M(I:y_1)=J\M(I:\M)\subseteq JI$. 
\end{proof}

 \section {The classes $\C$ and $\D$}
 
 In this section we define and study the properties of a class of $\M$-primary ideals of $R=K[x_1,\dots,x_n]$ denoted by $\C$ and of its subclass $\D$. Before giving the formal definition let us recall few notions that are needed in the sequel. 
Given an ideal $I$ with $\dim R/I=t$, the multiplicity $e(R/I)$ of $R/I$ is, by definition, $(t-1)!$ times the leading coefficient of the Hilbert polynomial of $R/I$ if $t>0$ and it is $\dim_K R/I$ otherwise. In particular, by definition, we have $e(R/R)=0$. 
 
 \begin{definition} We define $\C$ to be the class of the ideals $I$ of $R$ of finite colength such that:
\begin{itemize}
\item[(1)] $I + (\ell) = \M^{o(I)} + (\ell)$ for some non-zero linear form $\ell$,
\item[(2)] $I$ verifies one of the equivalent conditions of \ref{teoremone}.
\end{itemize}
We also set 
$$\D= \{ I \in \C \ : I \ {\mbox{ is integrally closed}}\}.$$ 
\end{definition}

\begin{remark}
\begin{itemize} 
\item[(1)] In the definition above we say ``finite colength" and not simply ``$\M$-primary" because we want $\C$ to contain $R$. 
\item[(2)] If $n=2$, then $\C$ is the class of contracted ideals.
\item[(3)] It follows from \cite[3.4]{CHH} that $\C$ can be also defined as the class of finite colength ideals $I$ which are componentwise linear with $\mu(I)= \mu(\M^{o(I)})$. 
\end{itemize} 
\end{remark} 

The next example shows that $\C$ cannot be defined as the class of ``contracted ideals with $ \mu(I)= \mu(\M^{o(I)})$". 

\begin{example} \label{bad} In $K[x_1,x_2,x_3]$ the ideal $I=(x_1^2,x_2x_3)+\M^3 $ is integrally closed, hence contracted and $\M$-full. Furthermore $\mu(I)= \mu(\M^{2})$. But $I \not \in \C$. 
\end{example}

However we have: 

\begin{lemma}\label{IIM}
Let $I$ be an ideal of $R$ of finite colength and order $d$. If both $I$ and $ \M I$ are $\M$-full and $ \mu(I)= \mu(\M^{d})$, then $\M I \in \C$. 
\end{lemma} 
\begin{proof} Since $I$ is $\M$-full and $ \mu(I)= \mu(\M^{d})$, then by \ref{mu} applied with $J=\M^{d}$ we deduce that there exists $\ell$ such that $\M^{d+1} + (\ell) = I \M + (\ell)$. 
 Since $\M I $ is $\M$-full, we conclude that $\M I\in \C$. 
\end{proof} 

 The class $\C$ is closed under the product and the integral closure. 
 \begin{proposition}
 \label{prodotto}
 If $I,J\in \C$, then $IJ\in C$ and $\overline I\in \D$. 
 \end{proposition}
\begin{proof} Set $d=o(I)$ and $d_1=o(J)$. Choose $\ell $ such that $I+(\ell)=\M^{d} +(\ell)$ and $J+(\ell)=\M^{d_1} +(\ell)$. Hence
 $\M^{d+d_1}+(\ell)\subseteq IJ+(\ell)$. Since the opposite inclusion is obvious, one has $\M^{d+d_1}+(\ell)= IJ+(\ell)$. Furthermore by \ref{elle} the dimension of $R/I_{\langle d\rangle}$ is $\leq 1$. Hence $IJ$ is componentwise linear by \ref{Id}. Hence $IJ\in \C$. As for $\overline I$ one notes that, by degree reasons, $o(\overline I)=d$ and $\M^d\subseteq I+(\ell)\subseteq \overline I+(\ell)$. Being integrally closed, $\overline I$ is contracted. It follows that $\overline I\in \D$. 
 \end{proof}
 
 Example \ref{exe2} shows that the class defined by the conditions ``$\M$-full and $\mu(I)= \mu(\M^{o(I)})$'', which properly contains $\C, $ is not closed under the product. 
 
 In dimension $2$, to every contracted ideal $I$ of order $d$ one associates its characteristic form $g(I)$ which is, by definition, the $\GCD$ of the elements in $I_{d}$. Zariski proved \cite[Appendix 5]{ZS} a factorization property for contracted ideals in dimension $2$. The factors are characterized by having pairwise coprime characteristic forms which are powers of irreducible forms. 
 Now we want to generalize Zariski's theorem to the class $\C$. To this end we will give another description of the ideals in it.

\begin{definition}
\label{q-family} We denote by $\A$ the set of the families $\Q=\{ Q_j\}_{j\in {\NN}}$ of homogeneous ideals of $R$ satisfying the following conditions: 
\begin{itemize} 
\item[(1)] $Q_j\subseteq Q_{j+1}$ for every $j$, 
\item[(2)] $Q_j=R$ for $j\gg 0$, 
\item[(3)] whenever $Q_j\neq R$, the ideal $Q_j$ is saturated and   $\dim R/Q_{j}=1$. 
\end{itemize} 
\end{definition} 
 
Given $\Q= \{ Q_i\}\in \A$, let $d_0=\reg(Q_0)$. For every $k\in \NN$ we set 
$$ I(\Q, k)= \oplus_{j \in \NN} (Q_j)_{d_0+k+j}. $$
We have: 
\begin{proposition}
\label{corrispondenza1}
For every $\Q=\{ Q_j\}\in \A$ and for every $k\in \NN$, one has $$ I(\Q, k) \in \C.$$
 \end{proposition}
\begin{proof}
Since $Q_j\subseteq Q_{j+1}$ we have $R_1Q_j\subseteq Q_{j+1}$ and hence $R_1(Q_j)_{d_0+j+k}\subseteq (Q_{j+1})_{d_0+j+1+k}$. This proves that $ I(\Q, k)$ is an ideal. If $Q_0=R$, then $I(\Q, k)=\M^k$ for all $k\geq 0$. 
Assume now that $Q_0\neq R$. Let $\ell$ be a linear form non-zero-divisor on $R/Q_0$. Since $\reg Q_0=d_0$, the ideal $Q_0+(\ell)$ is $0$-dimensional of regularity $d_0$. It follows that $\M^{d_0}\subseteq Q_0+(\ell)$. Therefore 
$\M^{d_0+k}\subseteq (Q_0)_{\langle d_0+k\rangle}+(\ell)\subseteq I(\Q, k)+(\ell)$ and hence
$\M^{d_0+k}+(\ell) = I(\Q, k)+(\ell)$. 
 It remains to prove that $ I(\Q, k)$ is componentwise linear, that is, $(Q_j)_{\langle d_0+h\rangle}$ has a linear resolution for every $h\in \NN$. By assumption $Q_j\subseteq Q_{j+1}$ and they define Cohen-Macaulay rings of the same dimension or are equal to $R$. It follows that $\reg(Q_j)\geq \reg(Q_{j+1})$. Hence $\reg(Q_j)\leq \reg(Q_0)= d_0$ for every $j$. Then for every $h\ge 0$ we have $(Q_j)_{\langle d_0+h\rangle}$ has a linear resolution. This proves the assertion. 
\end{proof}

\noindent Given an ideal $I$ in $ \C$ of order $d, $ for every  $j\geq 0, $ we set
$$Q_j(I)=(I_{\langle d+j\rangle})^{\sat}. $$ 

\begin{proposition}
\label{corrispondenza2}
Let $I\in \C$ and $d=o(I)$. For every $j\in \NN$ set $\Q(I)=\{Q_j(I)\}$ 
and $d_0=\reg(Q_0(I))$. 
Then $\Q(I)\in \A$ and $d\geq d_0$. 
\end{proposition} 

\begin{proof} 
Since $I_{\langle d\rangle}$ has dimension $\leq 1$, then $Q_j(I)$ is saturated of dimension $1$ or it is equal to $R$. Moreover 
$I_{\langle d+j\rangle}R_1\subseteq I_{\langle d+j+1\rangle}\subseteq I_{\langle d+j+1\rangle}^{\sat}$. Hence $Q_j(I)\subseteq Q_{j+1}(I)$. 
We have $d_0=\reg(Q_0(I))\leq \reg I_{\langle d
\rangle}=d$. 
\end{proof}

As a consequence we have:
\begin{theorem} 
\label{corrispondenza3}
With the notation of \ref{corrispondenza1} and \ref{corrispondenza2} the applications 
$$\varphi: \A\times \NN \longrightarrow \C \mbox{ and } \psi: \C \to \A \times \NN$$ 
 defined by $\varphi(\Q,k)=I(\Q,k)$ and $\psi(I)=(\Q(I), d-d_0)$ are inverse to each other. 
\end{theorem}
\begin{proof}
That the maps are well-defined follows from \ref{corrispondenza1} and \ref{corrispondenza2}. 
That are inverse to each other is a straightforward verification based on the observation that if $J$ is a saturated ideal generated in degree $\leq t, $ then $J_{\langle t\rangle}=J\cap \M^t$ and hence $J_{\langle t\rangle}^{\sat}=J$. 
\end{proof}

We need to recall now few facts about the ideal transform. Let $S=R[\M/\ell]$  where $\ell $ is a non-zero linear form. 
Clearly $\M S=(\ell)S$ and for every homogeneous element $f$ of degree $d$ one has $f=(f/\ell^d) \ell^d$ in $S$. Hence for every ideal $I$ of order $d$ we have 
$$IS=\ell^dI',$$ where $I'$ is an ideal of $S$. The ideal $I'$ is called the { ideal transform} of $I$ in
$S$.

\begin{proposition} 
\label{transform}
Let $I, J \in \C$ with $o(I)\geq o(J)$. The following facts are equivalent: 
\begin{itemize}
\item[(1)] $Q_j(I)=Q_j(J)$ for every $j$. 
\item[(2)] $I\M^s=J \M^r $ for some $r,s\in \NN$. 
\item[(3)] $I=J \M^r$ where $r=o(I)-o(J)$. 
\item[(4)] $I'=J'$ in $S=R[\M/\ell]$ for every linear form $\ell$. 
\item[(5)] $I'=J'$ in $S=R[\M/\ell]$ for a linear form $\ell$ not in $U(I)\cup U(J)$.
\end{itemize}
\end{proposition}
\begin{proof}
Conditions (1), (2) and (3) are equivalent by \ref{corrispondenza1}, \ref{corrispondenza2} and \ref{corrispondenza3}.
 That (3) implies (4) is clear by construction. That (4) implies (5) is obvious. 
Assume (5) and set $r=o(I)-o(J)$. Then $IS=\ell^{o(I)} I'$ and $J\M^{r}S=\ell^r\ell^{o(J)} J'=IS$. Since $J \in \C$, we have $J\M^{r}\in \C$ by \ref{prodotto}. Hence $I$ and $J\M^r$ are contracted from $S$. Since they have the same extension, it follows that $I=J\M^r$.
\end{proof}

\begin{definition} 
For $I,J\in \C$ we set $I\equiv J$ if $I$ and $J$ verify the equivalent conditions of \ref{transform}. 
\end{definition} 
In a different setting a similar equivalent relation is introduced  in \cite{L}. 

The extension $R\to R[\M/x_n]$ can be identified with the $K$-algebra homomorphism $\phi: R\to R$ sending $x_i\to x_ix_n$ for $i=1,\dots,n-1$ and $x_n$ to $x_n$. 
One has $\phi(f(x_1,\dots,x_n))=x_n^d f(x_1,\dots,x_{n-1},1)$ for every form of degree $d$. Denote by $\phi': R\to K[x_1,\dots,x_{n-1}]$ the dehomogenization map, that is, the $K$-algebra homomorphism sending $x_i\to x_i$ for $i=1,\dots,n-1$ and $x_n$ to $1$. So we have $\phi(f)=x_n^{d}\phi'(f)$ for every form of degree $d$. 

Let $I\in \C$ of order $d$. Let $P_1,\dots, P_m$ be the minimal primes of $Q_0(I)=I_{\langle d\rangle}^{\sat} $, necessarily homogeneous of dimension $1$ (with $m=0$ if $Q_0(I)=R$, that is, $I=\M^d$).
Note that, by construction, $I$ is contracted from any extension $R[\M/\ell]$ with $\ell\not\in \cup P_i$. 
After a change of coordinates, we may assume that $x_n \not\in \cup_{i=1}^m P_i$ and take $\ell=x_n$. We may write $I=\sum_{j\geq 0} I_{\langle j+d\rangle}$ and so 
$$\phi(I)R=\sum_{j\geq 0} \phi(I_{\langle j+d\rangle})R=x_n^d \sum_{j\geq 0} \phi'(I_{\langle j+d\rangle})x_n^j.$$ 
It follows that 
$$I'=\sum_{j\geq 0} \phi'(I_{\langle j+d\rangle}) x_n^j$$
that is 
$$I'= {\Big \{} \sum_j a_j x_n^j \ \ :\ \ a_j\in \phi'(I_{\langle j+d\rangle}){\Big \}} .$$
 
\begin{proposition} 
\label{maxmax}
With the notation above, we have: 
$$\sqrt{I'}=\cap_{i=1}^m (\phi'(P_i)R+(x_n))$$
and $\phi'(P_i)R+(x_n)$ are distinct maximal ideals of $R$. 
\end{proposition} 

\begin{proof} 
By definition, $Q_j(I)=I_{\langle j+d\rangle}^{\sat}$. Hence for some $u\in \NN$ one has $x_n^uQ_j(I)\subseteq I_{\langle j+d\rangle} \subseteq Q_j(I)$ which implies 
$$\phi'(I_{\langle j+d\rangle})=\phi'(Q_j(I)).$$
It follows that 
\begin{equation}
I'=\sum_{j\geq 0} \phi'(Q_j(I))x_n^j. \label{satdeo}
\end{equation} 
Since $Q_j(I)=R$ for $j\gg 0$ we have that $x_n^j\in I'$ for $j\gg 0$. As a consequence we have: 
$$\sqrt{I'}=\sqrt{\phi'(Q_0(I))R+(x_n)}=\sqrt{\phi'(Q_0(I))} R+(x_n)$$
 The known properties of the dehomogenization, see for instance \cite[Section 4.3]{KR}, guarantee that $\sqrt{\phi'(Q_0(I))}=\cap_{i=1}^m \phi'(P_i)$. The rest follows since $\phi'(P_i)$, as an ideal of $K[x_1,\dots,x_{n-1}]$, is maximal and $\phi'(P_i)\neq \phi'(P_j)$ for $i\neq j$. 
 \end{proof} 

The next result generalizes Zariski's factorization theorem for contracted ideals \cite[Appendix 5, Thm. 1]{ZS} to the class $\C$. The role played in \cite{ZS} by the characteristic form is played here by the ideal $Q_0(I)$. We call $Q_0(I)$ the {\em characteristic ideal} of $I$.

\begin{theorem}
\label{factorization}
Let $I\in \C$ and let $P_1, \dots, P_m$ be the minimal prime ideals of $Q_0(I)$. We have:
\begin{itemize} 
\item[(1)] There exist $L_1,\dots,L_m \in \C $ such that
$$ I \equiv L_1 L_2 \cdots L_m$$
and every $L_i$ has a $P_i$-primary characteristic ideal. 
\item[(2)] The $L_i$'s satisfying (1) are uniquely determined by $I$ up to $\equiv$. In particular, $Q_j(L_i)=Q_j(I)R_{P_i} \cap R$. 
\end{itemize} 
\end{theorem}
\begin{proof}
 First we prove that the $L_i$'s defined as in (2) satisfy (1) and then we prove the uniqueness of the $L_i$. For $i=1,\dots,m$ and $j\in \NN$ set $Q_i=\{Q_j(I)R_{P_i} \cap R \}_{j\in \NN}$. Then set $L_i=I(Q_i,0)$. By construction, $L_i\in \C$ and $Q_j(L_i)=Q_j(I)R_{P_i} \cap R$ and hence $Q_0(L_i)$ is $P_i$-primary. By \ref{prodotto} we have $L_1 L_2 \cdots L_m\in \C$. According to \ref{transform}, to prove (1) it is enough to show that 
 \begin{equation} 
 I'=L_1'L_2' \cdots L_m'
\label{pasta} 
\end{equation} 
in $S=R[\M/\ell] $ for a general linear form $\ell$. After a change of coordinates, we may assume that $x_n\not\in P_i$ for every $i$ and hence take $\ell=x_n$. Using formula (\ref{satdeo}) to describe $I'$ and the $L_i'$'s, (\ref{pasta}) becomes equivalent to 
 \begin{equation} 
\phi'(Q_j(I))=\sum_{*} \prod_{k=1}^m \phi'(Q_{j_k}(L_k))
\label{riso} 
\end{equation}
for all $j$, where the sum $\sum_{*}$ of the right hand side is extended to all the $j_1,\dots,j_m$ such that 
$j_1+j_2+\dots+j_m=j$. Equivalently, 
\begin{equation} 
\phi'(Q_j(I))= \phi' (\sum_* \prod_{k=1}^m Q_{j_k}(L_k)).
\label{minestra} 
\end{equation}
If we show that: 
\begin{claim} 
$Q_j(I)$ is the saturation of $\sum_* \prod_{k=1}^m Q_{j_k}(L_k)$
\end{claim} 
\noindent then we are done because two homogeneous ideals with the same saturation become equal after dehomogenization. 
 To prove the claim we localize $\sum_* \prod_{k=1}^m Q_{j_k}(L_k)$ at each $P_i$. What we get is $( \sum Q_{j_i}(L_i))R_{P_i}$ where the sum is exteded to $j_i\leq j$, that is, $Q_{j}(L_i)R_{P_i}$. Since $Q_j(L_i)=Q_j(I)R_{P_i} \cap R$ we have $Q_j(L_i)R_{P_i}=Q_j(I)R_{P_i}$. This proves the claim. 
 Now assume that there are other ideals $W_i\in C$ such that $I\equiv W_1\cdots W_m$ and $Q_0(W_i)$ is $P_i$-primary. Then $I'=W_1'\cdots W_m'$. Since by Proposition \ref{maxmax} the $W_i'$ are primary to distinct maximal ideals, we have that $I'=W_1'\cap \dots \cap W_m$ is a primary decomposition. By the uniqueness of minimal components in primary decompositions, we have $W_i'=L_i'$ and hence $W_i\equiv L_i$ as desired. 
 \end{proof}
 
We present now a formula for the Hilbert series of $I$ in terms of the Hilbert series of the ideals $ L_1,\dots, L_m$ appearing in the factorization of Theorem \ref{factorization}. If $\dim R=2, $ this has been already done in \cite[3.10]{CDJR}.
 
Since $I$ is an $\M$-primary ideal, then $\length(I^k/I^{k+1})$ is finite for every integer $k$. The Hilbert function $\HF_I(k)$ of $I$ is defined as
$$\HF_I(k)= \length(I^k/I^{k+1} ). $$ The Hilbert series of $I$ is 
$$
\HS_I(z)= \sum_{k\ge 0} \HF_I(k)z^k.
$$ 
It is well known that the Hilbert series is of the form 
$$\HS_I(z) = \frac {h_0(I)+h_1(I)z+\ldots+h_s(I)z^s}{(1-z)^n},$$ with $h_i(I)\in\ZZ$ for every $i$, $h_0(I)=\length(R/I)$ and $e(I)=\sum_{i=0}^s h_i(I)$ is the multiplicity of $I$. By definition, the h-polynomial of $I$ is $$h_I(z)=h_0(I)+h_1(I)z+\ldots+h_s(I)z^s.$$

\begin{lemma} 
\label{lunperI} Let $I$ be in $\C $ and let $I \equiv { L}_1 { L}_2 \cdots { {L}_m } $ be the factorization of \ref{factorization}. One has 
$$\length(\M^d/I)=\sum_{i=1}^m \length(\M^{d_i}/L_i) $$
where $d=o(I)$ and $d_i=o(L_i)$ for every $i=1,\ldots,m$. 
\end{lemma} 
\begin{proof} 
Since $\reg Q_j(I) \leq d$, then $\dim_K (R_{d+j}/I_{d+j})$ coincides with the multiplicity of $R/Q_j(I) $.
Hence $$\length (R/I)= \length (R/\M^d)+\sum_{i\ge 0} e(R/Q_j(I)).$$ 
Thus $\length (\M^d/I)=\sum_{j\ge 0} e (R/Q_j(I) )$. Since we know that 
$Q_j(I)=Q_j(L_1)\cap \dots \cap Q_j(L_m)$, the multiplicity formula \cite[4.7.8]{BH} implies that 
$e(R/Q_j(I))=\sum_{i=1}^m e(R/Q_j(L_i))$ and thus
\begin{align*}
\length (\M^d/I)={}&\sum_{j\ge 0} e(R/Q_j(I) )=\sum_{j\ge 0}\sum_{i=1}^m e(R/Q_j(L_i))\\ \\ 
={}&\sum_{i=1}^m \sum_{j\ge 0} e(R/Q_j(L_i)) =\sum_{r=1}^m \length (\M^{d_r}/L_r).
\end{align*}
\end{proof} 
 
\begin{proposition} 
\label{hfperI} 
With the notations of \ref{lunperI} we have: 
$$\HS_I(z)=\sum_{j=1}^m \HS_{L_j}(z)+\HS_{\M^d}(z)- \sum_{j=1}^m \HS_{\M^{d_j}}(z)
$$ and in particular 
$$e(I)=\sum_{j=1}^m e(L_j)+d^n - \sum_{j=1}^m d_j^n.$$ 
\end{proposition}
\begin{proof} 
 Note that for every integer $k$ the factorization of $I^k$ is: 
$$ I^k \ \equiv \ {L}_1^k \ { L}_2^k \cdots {L}_m^k $$ 
and hence 
$$\length(\M^{kd}/I^k)=\sum_{i=1}^m \length(\M^{kd_i}/L_i^k).$$
To conclude, first rewrite $\length(\M^{kd}/I^k)$ as $\length(R/I^k)-\length(R/\M^{kd})$ and similarly for the $L_i$'s and then sum up. 
 \end{proof} 

\begin{example}\label{3111} In $K[x,y,z]$ consider the ideal $I=(x^3,y^3,z^3,xy,yz,xz)$ of $\C$. We have $Q_0(I)=(xy,yz,xz)$ and $Q_j(I)=R$ for $j>0$. It follows from \ref{factorization} that $I\equiv L_1L_2L_3$ where $L_1=(x^2,y,z), L_2=(x,y^2,z), L_3=(x,y,z^2)$. To get an equality of ideals, we have to multiply the left hand side by $(x,y,z)$: 
$$ (x,y,z) (x^3,y^3,z^3,xy,yz,xz) = (x^2,y,z) (x,y^2,z) (x,y,z^2).$$
Taking into account that $d=2, d_1=d_2=d_3=1$ and that the $L_i$'s are complete intersections, we may apply \ref{hfperI} and get: 
$$\HS_I(z)=3\frac{2}{(1-z)^3}+\frac{ 4+4z}{(1-z)^3}-3\frac{1}{(1-z)^3}
$$
that is
$$\HS_I(z)=\frac{7+4z}{(1-z)^3}
$$
\end{example} 
\noindent The ideal of Example \ref{3111} appears in \cite{C} and \cite{L}. 
 \begin{theorem}
 \label{A}
Let $I\in \C$. Then 
\begin{itemize} 
\item[(1)] $\overline{\M I}=\M \overline{I}$. 
\item[(2)] $I\in \D$ if and only if $\M I\in \D$. 
\end{itemize} 
\end{theorem}
\begin{proof} (1) The inclusion $\M \overline{I}\subseteq \overline{\M I}$ holds in general, see \cite[1.1.3]{HS}. Using the characterization of integral closure by means of valuations, one shows that $$\overline {\M I}:\ell = \overline I $$ for every ideal $I$ and general linear form $\ell$, see the proof of \cite[3.1,3.3]{H2}  for details. Since $I \in \C$, then $ \M I +(\ell) = \M^{o(I)+1} + (\ell)$ is integrally closed. Then 
$\overline {\M I } \subseteq \overline {\M I + (\ell)}= {\M I + (\ell)}$. Hence 
$$ \overline {\M I}= ({\M I + (\ell)})\cap \overline {\M I}= \M I + \ell ( \overline {\M I} \ : \ \ell )= \M I+ \ell \overline I\subseteq \M \overline{I}.$$

(2) If $I\in \D$ then (1) implies $\M I\in \D$. Conversely if $\M I\in \D$ then $\overline {\M I}:\ell = {\M I}:\ell= \overline I$. Since $I$ is $\M$-full, it follows $I= \overline I.$
\end{proof}

Special cases of Theorem \ref{A} and Proposition \ref{ictransform} appear in \cite{DC1}. In general, even for a normal ideal $I$ the product $\M I$ need not be integrally closed, see \cite[Example 7.1]{DC2}. 

\begin{proposition} 
 \label{ictransform}
We have:
 \begin{itemize} 
 \item[(1)] If $I\in\D$ then $I'$ is integrally closed.
 \item[(2)] If $I'$ is integrally closed and $I$ is contracted, then $I$ is integrally closed.
 \end{itemize} 
 In particular if $I \in \C$, then $I\in \D$ if and only if $I'$ is integrally closed.
 \end{proposition}
 \begin{proof} Since $IS= \ell^d I' $ and $S$ is a polynomial ring (hence normal), then (1) follows if we prove that $IS$ is integrally closed. Consider the integral equation 
 $$ s^m + a_1 s^{m-1}+ \dots + a_m=0$$
 with $s \in S, $ $ a_i \in (I S)^i$. 
For every $i=0, \dots, m, $ we may write 
 $a_i = b_i /\ell^{\alpha} $ with $ b_i \in I^i \M^{\alpha} $ and $\alpha $ a fixed positive integer.
 Multiplying  by $\ell^{m\alpha} $ we get an equation among elements of $R$, namely
 $$ t^m + b_1 t^{m-1}+ \dots + (b_2 \ell^{\alpha}) t^{m-2} + \dots + (b_m/\ell^{\alpha}) =0$$
 where $t = s \ell^{\alpha} $ and $ b_i \ell^{(i-1)\alpha} \in I^i \M^{ i\alpha}$. Since $I \M^{\alpha} $ is integrally closed
by \ref{A}, it follows that $t=s\ell^{\alpha}\in I\M^{\alpha}$. Hence $s \in IS$. 

We prove now (2). Let $x \in R $ and $a_i \in I^i $  such that $$x^m +a_1 x^{m-1} + \dots +a_m=0 $$
and we claim that $x \in I$. Note that $a_i/\ell^{i d} \in (I')^i $ and 
$$ (x/\ell^d)^m + a_1/\ell^d (x/\ell^d)^{m-1} + \dots + a_m/\ell^{dm}=0. $$
Since $I'$ is integrally closed, it follows that $ x/\ell^d \in I'$, that is, $x\in IS$. Since $I$ is contracted we have $x\in I$. 
\end{proof}
 
\begin{theorem} 
\label{icfactorization}
Given  $I\in \C $ let $I\equiv L_1L_2 \cdots L_m $ be the factorization of \ref{factorization}. Then $I\in \D$ if and only if $L_j\in \D$ for every $j=1, \dots,m$.
\end{theorem}
\begin{proof} 
Assume that $I$ is integrally closed. By \ref{ictransform}(1), $I'=L_1'\cdots L_m'$ is integrally closed. Since $L_1',\dots, L_m'$ are primary to distinct maximal ideals, by localizing and contracting back one has that each $L_i'$ is integrally closed. Hence each $L_i$ is integrally closed by \ref{ictransform}(2). 
Conversely if $L_i$ is integrally closed for every $i=1, \dots, m$, then $L_i'$ is integrally closed by \ref{ictransform}(1). It follows that so is $I'=L_1'\cdots L_m'$ since $L_1'\cdots L_m'=L_1'\cap \cdots \cap L_m'$. Finally by \ref{ictransform}(2) we conclude that $I$ is integrally closed. 
\end{proof}

The following examples show that the class $ \D$ is not closed under product (for $n\geq 3$) and powers (for $n\geq 4$): 

\begin{example} \label{IJ} 
\begin{itemize} 
\item[(1)] The ideals $(x,y)^3+(x^2z)+\M^4$ and $(x,y)^3+(y^2z)+\M^4 $ of $K[x,y,z]$ are in $\D$ but not their product. 
\item[(2)] The ideal $(x^2,y^3,z^7,xy^2, xyz^2, xz^4,yz^5,y^2z^3,yz^5) \cap \M^7 +\M^8 $ of $K[x,y,z,t]$ is in $\D$ but not its square. 
\end{itemize}
\end{example}
 
Nevertheless, as an immediate consequence of Theorem \ref{icfactorization}, we have: 
 \begin{corollary} 
 \label{Mindistinct} 
 Let $I,J\in \D$ such that $Q_0(I)+Q_0(J)$ is $\M$-primary. Then $IJ\in \D$. 
 \end{corollary}
 
 Another corollary is: 
 
 \begin{corollary} 
 \label{cint-fact} 
 With the notation of \ref{factorization} we have 
 $$\overline{I} \equiv \overline{L_1}\ \overline{L_2} \cdots \overline{L_m}$$
 and $Q_0(\overline{L_i})$ is $P_i$-primary. 
 \end{corollary}
 \begin{proof} Combining \cite[Exercise 1.1, p. 20]{HS} with \ref{A}, we get: 
 $$\overline{I} \equiv \overline{ \overline{L_1}\ \overline{L_2} \cdots \overline{L_m} }.$$
 The conclusion follows from \ref{icfactorization} provided we prove that $Q_0(\overline{L_i})$ is $P_i$-primary. So assume that $L\in \C$ has order $d$ and $Q_0(L)$ is $P$-primary for some $1$-dimensional prime $P$. Set $J=\overline{L}$. By degree reasons, $J_d \subset \overline{ L_{\langle d \rangle}}$ and 
 $\overline{L_{\langle d\rangle}}\subseteq \sqrt{L_{\langle d\rangle}}=P$. Hence 
 $J_{\langle d\rangle}\subseteq P$ which implies that $Q_0(J)$ is $P$-primary. 
\end{proof}

 \section{The Goto-classes $\G$ and $\G^*$}
\label{casilimite}
Consider the following subclass of $\C$: 
\begin{definition} We define the Goto-class $\G$ to be the set of the ideals $I\in \C$ such that:
\begin{itemize} 
\item[(1)] The minimal primes $P_1,\dots,P_m$ of $Q_0(I)$ are geometrically prime, equivalently, each $P_i$ is generated by $n-1$ linearly independent linear forms in $R=K[x_1, \dots, x_n]$ (e.g. $K$ is algebraically closed). 
\item[(2)] For every $j\in \NN$ the primary components of $Q_j(I)$ are powers of the $P_i$'s. That is, 
$$Q_j(I)=\cap_{i=1}^m {\ {P}}_i^{\alpha_{ij }}$$
with $\alpha_{ij}\in \NN$. 
\end{itemize} 
Further we set: 
$$\G^*=\{ I\in \G : I \mbox{ is integrally closed }\}.$$
\end{definition}

In dimension two $\G=\C$ and it coincides with the whole class of contracted ideals. Our goal is to show that the Goto-classes $\G$ and $\G^*$ behave, to a certain extent and respectively, as the class of contracted ideals and the class of integrally closed ideals in dimension $2$. The factorization in Theorem  \ref{factorization} will allow to reduce most of the problems to the case of ideals in $\G$ with a primary characteristic ideal. So we will discuss in some details the properties of these ideals. 

Let $P$ be a geometrically prime ideal of $R$ of dimension $1$. Let $L\in \G$ of order $d$ such that $Q_0(L)$ is $P$-primary. Then $Q_j(L)=P^{\alpha_j}$ where the $\alpha_j$'s form a weakly decreasing integral sequence with $\alpha_j=0$ for $j\gg 0$. Hence $L$ is described by the triplet $P,\{\alpha_j\}$ and $d$. We give another description of $L$ that best suits our needs. Indeed, one  shows that there exists a uniquely determined sequence of integers $0=a_0<a_1<\dots<a_d$ such that 
\begin{equation}
\label{eccolo}
L=\sum_{i=0}^d P^{d-i}\ell^{a_i}
\end{equation}
where $\ell$ is any linear form not in $P$. 
To emphasize the dependence of $L$ on $P$ and the sequence $a_0,\dots,a_d$ we will denote $L$ by $L(P,a)$, that is, 
\begin{equation}
\label{rieccolo}
L(P,a)=\sum_{i=0}^d P^{d-i}\ell^{a_i}
\end{equation}
\begin{example}
Let $R=K[x_1,x_2,x_3]$ and $P=(x_1,x_2)$. Associated with the sequence $\alpha=(5,3,3,2,0,0\dots)$ and with $d=6$ we have the ideal $L$ whose components are $L_{6+j}=(P^{\alpha_j})_{6+j}$ for $j\geq 0$. 
We can write $L$ as $L(P,a)=\sum P^{d-i}x_3^{a_i} $ where $a=(0,1,3,4,7,9,10)$. 
\end{example}

Given two sequences of integers $a=(a_0,\dots, a_d)$ and $b=(b_0,\dots,b_e)$ we define their product $ab$ to be the sequence $(c_0,\dots,c_{d+e})$ where  $c_j=\min\{ a_r+b_s : r+s=j\}$. Furthermore we denote by $a^{(k)}$ the product of $a$ with itself $k$ times. 
By the very definition one has: 
$$L(P,a)L(P,b)=L(P,ab) \mbox{ and } L(P,a)^k=L(P, a^{(k)})$$
for every $a,b$ and $P$. We have: 

\begin{proposition}
\label{DM270} 
Let $a=(a_0,\dots, a_d) \in \NN^{d+1}$ be an increasing sequence with $a_0=0$.
\begin{itemize} 
\item[(1)] There exists an increasing sequence $a'=(a'_0,\dots,a'_d)$ with $a'_0=0$ (uniquely determined by $a$) such that for every $n>1$ and for every $1$-dimensional geometrically prime ideal $P$ of $R=K[x_1,\dots,x_n]$ one has $\overline{L(P,a)}=L(P,a')$. 

\item[(2)] The following conditions are equivalent: 
\begin{itemize}
\item[(i)] $L(P,a)$ is integrally closed 
 for every $n>1$ and for every $1$-dimensional geometrically prime ideal $P$ of $R$.
 \item[(ii)] $L(P,a)$ is integrally closed for some $n>1$ and some $1$-dimensional geometrically prime ideal $P$ of $R$. 
 \item[(ii)] $a=a'$. 
\end{itemize}
\end{itemize}
\end{proposition}
\begin{proof} (1) Let $n>1$ and let $P$ be a $1$-dimensional geometrically prime ideal of $R$. Choosing bases properly, we may assume that $P=(x_1,\dots,x_{n-1})$ and $\ell=x_n$ so that $L(P,a)$ is a monomial ideal. The integral closure of a monomial ideal $I$ is the ideal generated by the monomials $m$ such that $m^k\in I^k$ for some $k>0$. A monomial $m=m_1x_n^{d-j}$ with $m_1$ supported on $x_1,\dots,x_{n-1}$ satisfies $m^k=m_1^k x_n^{kd-kj} \in L(P,a)^k$ iff $\deg m_1^k=k \deg m_1 \geq (a^{(k)})_{kj} $ iff $\deg m_1\geq (a^{(k)})_{kj}/k$. Hence setting $$a'_j=\min\{\lceil (a^{(k)})_{kj}/k \rceil :k>0\}$$ we get (1). Statement (2) follows immediately from (1). \end{proof} 

Given a $1$-dimensional geometrically prime ideal $P$ of $R$ and numbers $d,t\in \NN$ with $d\le t$ we set 
$$J_P(d,t)= \overline {P^d+\M^t},$$
equivalently 
$$J_P(d,t)= \overline {(\ell_1^d,\dots,\ell_{n-1}^d, \ell^t)},$$
where $P=(\ell_1,\dots,\ell_{n-1})$ and $\ell$ is a linear form not in $P$. 
By construction, $J_P(d,t)\in \G^*$ and its characteristic ideal is $P^d$ unless $t=d$. Hence $J_P(d,t)$ must be of the form $L(P,a)$ for a sequence $a$. 
Indeed a simple computation shows that: 
$$J_P(d,t)=L(P,a)$$
where $a=(a_0,\dots,a_d)$ with $a_i=\lceil it/d\rceil$ for $i=0,\dots,d$.

We say that an ideal $I$ is simple if it cannot be written as a product of proper ideals. 
\begin{remark}
\label{easy?}
It is an easy exercise and part of the folklore of the subject that $\overline{(x_1^d,x_2^t)}$ is simple in $K[x_1,x_2]$ iff $\GCD(d,t)=1$ and that every simple integrally closed ideal of $K[x_1,x_2]$ with characteristic form equal to $x_1$ is of the form $\overline{(x_1^d,x_2^t)}$. 
\end{remark}

\begin{proposition}
\label{DM271} 
Let $a=(a_0,\dots, a_d) \in \NN^{d+1}$ be an increasing sequence with $a_0=0$ and $P$ a $1$-dimensional geometrically prime ideal of $R$. Then the following conditions are equivalent: 
\begin{itemize}
\item[(1)] $L(P,a)$ is integrally closed, simple and different from $\M$. 
\item[(2)] there exists $t>d$ such that $\GCD(d,t)=1$ and $L(P,a)=J_P(d,t)$. 
\item[(3)] there exists $t>d$ such that $\GCD(d,t)=1$ and $a_i=\lceil it/d\rceil$ for $i=0,\dots,d$.
\end{itemize}
\end{proposition} 
\begin{proof} 
The result follows from \ref{DM270},  \ref{easy?} and  the following claim: 
\begin{claim} 
$L(P,a)$ is integrally closed and simple in $R= K[x_1, \dots, x_n] $ if and only if $L((x_1),a)$ is integrally closed and simple in $K[x_1,x_2]$. 
\end{claim}
To prove the claim assume first that $L(P,a)$ is integrally closed and simple. Then $L((x_1),a)$ is integrally closed by \ref{DM270}. If, by contradiction, $L((x_1),a)$ is not simple, then $L((x_1),a)=IJ$ with $I,J$ integrally closed. Hence $I$ and $J$ are 
of the form $I=L((x_1),b)$ and $J=L((x_1),c)$. 
It follows that $L(P,a)=L(P,b)L(P,c)$ contradicting the fact that $L(P,a)$ is simple. 

Viceversa, assume that $L((x_1),a)$ is integrally closed and simple. Then $L(P,a)$ is integrally closed by \ref{DM270}. If, by contradiction, $L(P,a)$ is not simple, then $L(P,a)=IJ$ with $I,J$ proper ideals. Since $\M^u\subset L(P,a)\subseteq I$ it follows that $\sqrt I=\M$ and for the same reason $\sqrt J=\M$. 
After a change of coordinates, we may assume that $P=(x_1,\dots,x_{n-1})$ and consider the $K$-algebra homomorphism $\psi: R\to K[x_1,x_2]$ sending $x_i$ to $x_1$ for $i<n$ and $x_n$ to $x_2$. We have $L((x_1),a)=\psi(L(P,a))=\psi(I)\psi(J)$ and $\psi(I)$ and $\psi(J)$ are proper since $\psi(I)\subseteq \psi(\M)=(x_1,x_2)$ and similarly for $\psi(J)$. This contradicts the assumptions and proves the claim.
\end{proof} 

Next we show that the factorization of Theorem \ref{factorization} restricts to the class $\G$.

\begin{proposition}
\label{1/2G}
We have:
\begin{itemize} 
\item[(1)] Let $I\in \C$ be such that the minimal primes of $Q_0(I)$ are geometrically prime. Let $I\equiv L_1\cdots L_m$ be the factorization of \ref{factorization}. Then $I\in \G$ iff $L_i\in \G$ for every $i$. 
\item[(2)] $\G$ is closed under product. 
\item[(3)] If $I\in \G$ then $\overline I\in \G^*$. 
\end{itemize} 
\end{proposition} 
\begin{proof} (1) Let $\{P_1,\dots,P_m\}$ the minimal primes of $Q_0(I)$. By \ref{factorization} we know that $Q_j(I)R_{P_i}=Q_j(L_i)R_{P_i}$ and this implies the assertion. 

(2) Let $I,J\in \G$. Set $d=o(I)$ and $c=o(J)$. We have to show that for every $P\in \Min(Q_0(IJ))$ and for every $j$ we have $(IJ)_{\langle d+c+j\rangle}R_P$ is a power of $PR_P$. 
Note that $(IJ)_{\langle d+c+j\rangle}=\sum_{i=0}^j I_{\langle d+i\rangle}J_{\langle c+j-i\rangle}$ and that $I_{d+i}R_P=P^{a_i}R_P$ and $J_{\langle c+j-i\rangle}R_P=P^{b_{j-i}}R_P$ for non-negative integers $a_i$ and $b_i$. It follows that we have $(IJ)_{\langle d+c+j\rangle}R_P=P^tR_P$ where $t=\min\{ a_i+b_{j-i} : i=0,\dots,j\}$. 

(3) Let $I\equiv L_1\cdots L_m$ be the factorization of \ref{factorization}. Then by \ref{cint-fact} 
 we have $\overline{I} \equiv \prod_i \overline{L_i}$. By (1) it is enough to show that $\overline{L_i}\in \G$. We may hence assume that $I$ is of the form $L(P,a)$. But we have already observed in \ref{DM270} that $\overline{L(P,a)}=L(P,a')$, which implies that $\overline{L(P,a)}\in \G$. 
\end{proof} 

We can state now the main result of the section.

\begin{theorem}
\label{mainG}
 We have:
 \begin{itemize}
 \item[(1)] $\G^*$ is closed under product. In particular, every $I\in \G^*$ is normal. 
\item[(2)] Every $I\in \G^*$ has a factorization $$I\equiv J_1\cdots J_t$$ where $J_i\in\G^*$ is simple and $Q_0(J_i)$ is primary for every $i=1,\dots,t$. 
\item[(3)] In the factorization of (2), the factors $J_i$ are uniquely determined by $I$ up to order. Moreover, each $J_i$ is of the form $J_{P_i}(d_i,t_i)$ and $d_i<t_i$ with $\GCD(d_i,t_i)=1$. 
\end{itemize} 
\end{theorem} 
\begin{proof} 
(1) Let $I,J\in \G^*$. By \ref{1/2G}(2) we know that $IJ\in \G$ and we have to prove that $IJ$ is integrally closed. By \ref{factorization}, \ref{1/2G}(1) and \ref{icfactorization} we have factorizations $I\equiv L_1\cdots L_m$ and $J\equiv U_1\cdots U_r$ and the $L_i$ and $U_j$ belong to $\G^*$. 
Hence $IJ\equiv \prod L_i\prod U_i$. If $L$ and $U$ have $P$-primary characteristic ideal then the same is true for $LU$. Hence, the factors in the (unique) factorization of \ref{factorization} of $IJ$ are of the form $L_iU_j$ (if $L_i $ and $U_j$ have $ P$-primary characteristic ideal with respect to the same prime) or $L_i$ or $U_j$. By virtue of \ref{icfactorization} we may assume right away that $I$ and $J$ have $P$-primary characteristic ideal, say $I=L(P,a)$ and $J=L(P,b)$. Then $IJ=L(P,ab)$. Since $I$ and $J$ are integrally closed, the same is true for $L((x_1),a)$ and $L((x_1),b)$ in $K[x_1,x_2]$ by \ref{DM270}. As in dimension $2$ the product of integrally closed ideals is integrally closed, we have that $L((x_1),a)L((x_1),b)=L((x_1),ab)$ is integrally closed. By \ref{DM270} it follows that $IJ$ is integrally closed. 

(2) By virtue of \ref{1/2G} we have $I\equiv L_1\cdots L_m$ with $L_i\in \G^*$ and $Q_0(L_i)$ primary.
 Hence we may assume that $I=L(P,a)$ for some $1$-dimensional geometrically prime ideal $P$ and a sequence $a=(a_0,\dots,a_d)$. By Zarisky factorization theorem \cite{ZS} and  \ref{easy?}  one has $L((x_1),a)=(x_1,x_2)^c J_{(x_1)}(d_1,t_1)\cdots J_{(x_1)}(d_p,t_p)$ with $d_i<t_i$ and $\GCD(d_i,t_i)=1$. It follows that $I=\M^c J_{P}(d_1,t_1)\cdots J_{P}(d_p,t_p)$ and hence $I\equiv J_{P}(d_1,t_1)\cdots J_{P}(d_p,t_p)$. The conclusion follows from \ref{DM271}. 
 
 (3) That the factors of the factorization in (2) are of the form $J_{P_i}(d_i,t_i)$ with $d_i<t_i$ and $\GCD(d_i,t_i)=1$ has been already proved. It remains to prove the uniqueness. Suppose we have two factorizations of $I$ as in (2). By \ref{1/2G} and \ref{factorization} we may assume that the characteristic form of $I$ is $P$-primary. Hence we have $I\equiv J_P(d_1,t_1)\cdots J_P(d_p,t_p)$ and  $I\equiv J_P(c_1,s_1)\cdots J_P(c_q,s_q)$ with $d_i<t_i$ and $\GCD(d_i,t_i)=1$ as well as $c_i<s_i$ and $\GCD(c_i,s_i)=1$. As a consequence we have $\M^a J_P(d_1,t_1)\cdots J_P(d_p,t_p)=\M^b J_P(c_1,s_1)\cdots J_P(c_q,s_q)$, and it follows that $(x_1,x_2)^a J_{(x_1)}(d_1,t_1)\cdots J_{(x_1)}(d_p,t_p)=(x_1,x_2)^b J_{(x_1)}(c_1,s_1)\cdots J_{(x_1)}(c_q,s_q)$ in $K[x_1,x_2]$. By the uniqueness of the factorization of integrally closed ideals in $K[x_1,x_2]$, we have that $p=q$ and, up to the order, 
 $(d_i,t_i)=(c_i,s_i)$ for $i=1,\dots,p$. Hence $J_P(d_i,t_i)=J_P(c_i,s_i)$ for $i=1,\dots,p$ proving the assertion. 
\end{proof}

\begin{remark}\label{easy}
Let $I$ be an $\M$-primary complete intersection ideal of $R$. Goto proved in \cite{G} that the following conditions are equivalent: 
\begin{itemize} 
\item[(1)] $I$ is integrally closed.
\item[(2)] $I$ is normal.
\item[(3)] $I=(\ell_1,\dots,\ell_{n-1},\ell_n^t)$ for linearly independent linear forms $\ell_1,\dots,\ell_n$ and some $t>0$.
\end{itemize}
Complete intersections satisfying these equivalent conditions are called of Goto-type (see \cite{CHV}). 
Note that the ideals of Goto-type are in the Goto-class $\G$, they are exactly the ideals of type 
$J_P(1,t)$ used above. 
\end{remark} 
 
Hence, as a consequence of Theorem \ref{mainG}, we have: 
\begin{corollary} 
\label{goto-normal}
The product of complete intersections of Goto-type is a normal ideal.
\end{corollary}
 
In dimension $2$, every integrally closed ideal has a Cohen-Macaulay associated graded ring (see \cite{H2}, \cite{LT}). This is no longer true in higher dimension and not even for normal ideals. The first examples of  normal ideals with non Cohen-Macaulay associated graded ring is given by a construction of  Cutkosky  \cite{C3}. Later on Huckaba and Huneke \cite[Theorem 3.12]{HuHu}) proved that 
$$I=(x^4)+(x,y,z)(y^3+z^3)+(x,y,z)^5\subseteq K[x,y,z]$$ 
 is normal, but $\gr_{I^n}(R)$ is not Cohen-Macaulay for every $n$. 

One might, however, ask: 

\begin{question} 
\label{que}
Let $I\in \G^*$. Is $\gr_I(R)$ Cohen-Macaulay? 
\end{question}
We show that Question \ref{que} has a positive answer in two cases. The first is the following.
 \begin{corollary}
\label{megliocheniente}
 Let $I\in \G^*$. Then $\Rees(I)$ is normal. In particular, $\Rees(I)$, equivalently $\gr_I(R)$, is Cohen-Macaulay if $I$ is monomial in some system of coordinates (e.g.~the characteristic ideal of $I$ has at most $2$ minimal primes). 
\end{corollary}
\begin{proof} The first assertion follows from \ref{mainG}(1). The second follows from the fact that if 
the characteristic ideal of $I$ has at most $2$ minimal primes, then up to a choice of coordinates, we may assume that $I$ is monomial. For a monomial ideal $I$, the normality of $\Rees(I)$ implies its Cohen-Macaulayness as proved by Hochster \cite[6.3.5¤]{BH}. 
\end{proof} 
To show that \ref{que} has a positive answer if $\dim R \le 3$ we need the following result. 
 \begin{lemma} 
 \label{hhh}
 Let $I$ be an $\M$-primary ideal of $R=K[x_1, \dots, x_n]$. If $\gr_I(R)$
 is Cohen-Macaulay, then the degree of its $h$-polynomial is $\leq n-1$.
 \end{lemma}
 \begin{proof} 
Since the ideal $I$ is $\M$-primary, then $\gr_I(R)$ is Cohen-Macaulay if and only if $\gr_{I_{\M}}(R_{\M})$ is Cohen-Macaulay; moreover $\gr_I(R)$ and $\gr_{I_{\M}}(R_{\M})$ have the same Hilbert series. Hence we may reduce the problem to the local case (see for example Remark 2.2. \cite{CDJR}). Note that if $J$ is a minimal reduction of $I$, then the $h$-polynomial $h_I(z)=h_0(I)+h_1(I)z+\ldots+h_s(I)z^s $ coincides with the Hilbert series of the ideal $I/J$. Now by a consequence of Briancon-Skoda \cite[11.1.9]{HS}, we have $I^n\subseteq J, $ hence $\HF_{I/J}(n)= \length(I^n +J/I^{n+1}+J) =h_n(I)=0 $ and the result follows.
\end{proof}
 
 \begin{theorem} 
 \label{3CM}
 Assume that $\dim R\leq 3$. If $I\in \G^*$, then $\gr_I(R)$ is Cohen-Macaulay.
 \end{theorem}
 
 \begin{proof} Consider the factorization $I \equiv L_1 \cdots L_m$ of \ref{factorization}. We know by \ref{icfactorization} and \ref{1/2G} that $L_i\in \G^*$. By \ref{hfperI} one has $h_I(z)=\sum_{j=1}^m h_{L_j}(z)+h_{\M^d}(z)- \sum_{j=1}^m h_{\M^{d_j}}(z)$. By \ref{megliocheniente} we know that $\gr_{L_i}(R)$ is 
Cohen-Macaulay for every $i=1, \dots, m$. That $\gr_{\M^u}(R)$ is Cohen-Macaulay for every $u$ is well-know. Thus by \ref{hhh} the degree of $h_{L_i}(z) $ $\le 2$ for every $i=1,\dots,m$ and the same is true for $h_{\M^u}(z)$. It follows that the degree of $h_I(z)$ is $\le 2$.
Localizing at $\M$ we may assume that $R$ is local. Let $J$ be a minimal reduction of $I;$ since $I$ is integrally closed, by \cite{H1} we have $I^2 \cap J=JI$. Then the result follows by \cite[2.2]{GR}.
\end{proof}

\end{document}